\documentclass[twoside,leqno]{amsart}

\usepackage{amsmath,color,amsthm}

\numberwithin{equation}{section}

\newtheorem{theorem}{Theorem}[section]
\newtheorem{corollary}[theorem]{Corollary}
\newtheorem{lemma}[theorem]{Lemma}
\newtheorem{proposition}[theorem]{Proposition}

\newtheorem{example}[theorem]{\sl Example}
\newtheorem{definition}[theorem]{\sl Definition}

\theoremstyle{definition}
\newtheorem{Remark}[theorem]{Remark}

\newcommand{\beqn}{\begin{eqnarray}}
\newcommand{\eeqn}{\end{eqnarray}}
\newcommand{\beqnn}{\begin{eqnarray*}}
\newcommand{\eeqnn}{\end{eqnarray*}}

\newcommand{\EE}{{\bf  E}}

\newcommand{\ct}{{\tilde{c}}}

\newcommand{\Lc}{{\mathcal L}}

\newcommand{\Lpto}{{\stackrel{L^p}{\longrightarrow}}}

\newcommand{\begp}{\begin{proposition}}
\newcommand{\enp}{\end{proposition}}
\newcommand{\begt}{\begin{theorem}}
\newcommand{\ent}{\end{theorem}}
\newcommand{\begl}{\begin{lemma}}
\newcommand{\enl}{\end{lemma}}
\newcommand{\begc}{\begin{corollary}}
\newcommand{\enc}{\end{corollary}}
\newcommand{\begcl}{\begin{claim}}
\newcommand{\encl}{\end{claim}}
\newcommand{\begr}{\begin{Remark}\rm}
\newcommand{\enr}{\end{Remark}}
\newcommand{\begal}{\begin{algorithm}}
\newcommand{\enal}{\end{algorithm}}
\newcommand{\begd}{\begin{definition}}
\newcommand{\enf}{\end{definition}}
\newcommand{\begx}{\begin{example}}
\newcommand{\enx}{\end{example}}
\newcommand{\bega}{\begin{array}}
\newcommand{\ena}{\end{array}}

\def\rompar(#1){\textup(#1\textup)}    

\newcommand\ga{\alpha}

\newcommand{\refS}[1]{Section~\ref{#1}}
\newcommand{\refT}[1]{Theorem~\ref{#1}}

\newcommand{\refL}[1]{Lemma~\ref{#1}}
\newcommand{\refP}[1]{Proposition~\ref{#1}}

\newcommand\QuickSort{\texttt{QuickSort}}

\newcommand{\ignore}[1]{}

\begin{document}

\newcommand{\tab}[0]{\hspace{.1in}}

\title[Exact $L^2$-Distance for QuickSort Key Comparisons]
{Exact $L^2$-Distance from the Limit for QuickSort Key Comparisons (Extended Abstract)}

\author{\ \ \ \ Patrick Bindjeme\ \ \ \ \vspace{.2cm}\\ James Allen Fill}
\address{Department of Applied Mathematics and Statistics,
The Johns Hopkins University,
34th and 
Charles Streets,
Baltimore, MD 21218-2682 USA}
\email{bindjeme@ams.jhu.edu and jimfill@jhu.edu}
\thanks{Research supported by the Acheson~J.~Duncan Fund for the Advancement of Research in Statistics.}

\date{January~26, 2012}

\maketitle

\begin{center}
{\sc Abstract}
\vspace{.3cm}
\end{center}

\begin{small}
Using a recursive approach, we obtain a simple exact expression for the $L^2$-distance from the limit in 
R\'{e}gnier's~\cite{r1989} classical limit theorem for the number of key comparisons required by {\tt QuickSort}.  A previous study by Fill and Janson~\cite{fj2002} using a similar approach found that the $d_2$-distance is of order between $n^{-1} \log n$ and $n^{-1/2}$, and another by Neininger and Ruschendorf~\cite{neiru2002} found that the Zolotarev $\zeta_3$-distance is of exact order $n^{-1} \log n$.  Our expression reveals that the $L^2$-distance is asymptotically equivalent to $(2 n^{-1} \ln n)^{1/2}$.
\end{small}
\bigskip

\section{Introduction, review of related literature, and summary}\label{S:introduction}
We consider Hoare's~\cite{h1962} \QuickSort\ sorting algorithm applied to an infinite stream of iid (independent and identically distributed) uniform random variables $U_1, U_2, \dots$.  \QuickSort\ chooses the first key $U_1$ as the ``pivot'', compares each of the other keys to it, and then proceeds recursively to sort both the keys smaller than the pivot and those larger than it.  If, for example, the initial round of comparisons finds $U_2 < U_1$, then $U_2$ is used as the pivot in the recursive call to the algorithm that sorts the keys smaller than $U_1$ because it is the first element in the sequence $U_1, U_2, \dots$ which is smaller than $U_1$.  In a natural and obvious way, a realization (requiring infinite time) of the algorithm produces an infinite rooted binary search tree which with probability one has the \emph{completeness} property that each node has two child-nodes.  

Essentially the same algorithm can of course be applied to the truncated sequence $U_1, U_2, \dots, U_n$ for any finite~$n$, where the recursion ends by declaring that a list of size~$0$ or~$1$ is already sorted.  Let $K_n$ denote the number of key comparisons required by \QuickSort\ to sort $U_1, U_2, \dots, U_n$.  Then, with the way we have set things up, all the random variables $K_n$ are defined on a common probability space, and $K_n$ is nondecreasing in~$n$.  Indeed, $K_n - K_{n - 1}$ is simply the cost of inserting $U_n$ into the usual (finite) binary search tree formed from $U_1, \dots, U_{n - 1}$.

In this framework, R\'{e}gnier~\cite{r1989} used martingale techniques to establish the following $L^p$-limit theorem; she also proved almost sure
convergence.  We let
$$
\mu_n := \EE\,K_n.
$$

\begt[R\'{e}gnier \cite{r1989}]
\label{thregn}
There exists a random variable $T$ satisfying
\begin{eqnarray*}
Y_n := \frac{K_n - \mu_n}{n+1}\,\Lpto\,T
\end{eqnarray*}
for every finite $p$.
\ent

R\"{o}sler~\cite{r1991} characterized the distribution of R\'{e}gnier's limiting~$T$ as the unique fixed point of a certain distributional transformation, but he also described explicitly how to construct a random variable having the same distribution as~$T$.  We will describe his explicit construction in equivalent terms, but first we need two paragraphs of notation.

The nodes of the complete infinite binary search tree are labeled in the natural binary way: the root gets an empty label written $\varepsilon$ here, the left (respectively, right) child is labeled 0 (resp.,\ 1), the left child of node~0 is labeled~00, etc.  We write $\Theta := \cup_{0 \leq k < \infty} \{0, 1\}^k$ for the set of all such labels.  If $V_\theta$ denotes the key inserted at node $\theta \in \Theta$, let $L_\theta$ (resp.,\ $R_\theta$) denote the largest key smaller than $V_\theta$ (resp.,\ smallest key larger than $V_\theta$) inserted at any ancestor of $\theta$, with the exceptions $L_\theta := 0$ and $R_\theta := 1$ if the specified ancestor keys do not exist.
Further, for each node~$\theta$, define
\begin{align}
\phi_{\theta} &:= R_{\theta} - L_{\theta}, \qquad U_{\theta} := \phi_{\theta 0} / \phi_{\theta}, \nonumber \\ 
\label{gtheta}
G_{\theta} &:= \phi_{\theta} C(U_{\theta}) = \phi_{\theta} - 2\phi_{\theta} \ln \phi_{\theta} + 2\phi_{\theta 0} \ln \phi_{\theta 0} + 2\phi_{\theta 1} \ln \phi_{\theta 1},
\end{align}
where for $0 < x < 1$ we define
\begin{equation}
\label{Cdef}
C(x) := 1 + 2 x \ln x + 2 (1 - x) \ln(1 - x).
\end{equation}

Let $1 \leq p < \infty$.  The \emph{$d_p$-metric} is the metric on the space of all probability distributions with finite 
$p$th absolute moment defined by
\begin{eqnarray*}
d_p(F_1, F_2) & := & \inf \|X_1 - X_2\|_p,
\end{eqnarray*}
where we take the infimum of $L^p$-distances over all pairs of random variables $X_1$ and $X_2$ (defined on the same probability space) with respective marginal distributions $F_1$ and $F_2$.  By the $d_p$-distance between two random variables we mean the $d_p$-distance between their distributions. 

We are now prepared to state R\"{o}sler's main result.  {\sc Note}:\ Here and later results have been adjusted slightly as necessary to utilize the same denominator $n + 1$ (rather than $n$) that R\'{e}gnier used.

\begt[R\"{o}sler~\cite{r1991}]\label{throsl}
For any finite~$p$, the infinite series $Y =  \sum_{j=0}^\infty \sum_{|\theta|=j}G_\theta$ converges in~$L^p$, and the
sequence $Y_n = (K_n - \mu_n) / (n + 1)$ converges in the $d_p$-metric to~$Y$. 
\ent

Of course it follows from Theorems~\ref{thregn}--\ref{throsl} that~$T$ and~$Y$ have the same distribution.  {\bf The purpose of the present extended abstract is to show that in fact $T = Y$ and to provide a simple explicit expression for the $L^2$-distance between $Y_n$ and~$Y$} valid for every~$n$; this is done in \refT{L2ratek} below.

We are aware of only two previous studies of the rate of convergence of $Y_n$ to~$Y$, and both of those concern certain distances between \emph{distributions} rather than between \emph{random variables}.  The first study, by
Fill and Janson \cite{fj2002}, provides upper and lower bounds on $d_p(Y_n, Y)$ for general~$p$; we choose to focus here on $d_2$. 
\begt[Fill and Janson \cite{fj2002}]
\label{thfjd2} 
There is a constant $c > 0$ such that for any $n \geq 1$ we have
$$
c n^{-1} \ln n \le d_2(Y_n, Y) < 2n^{- 1 / 2}.
$$
\ent
\noindent
To our knowledge, the gap between the rates $(\log n) / n$ and $n^{-1/2}$ has not been narrowed.
Neininger and Ruschendorff~\cite{neiru2002} used the Zolotarev $\zeta_3$-metric and found that the correct rate in that metric is $n^{-1} \log n$, but their techniques are not sufficiently sharp to obtain $\zeta_3(Y_n, Y) \sim \ct n^{-1} \ln n$ for some constant~$\ct$.

In our {\bf main \refT{L2ratek}}, proved using the same recursive approach as in Fill and Janson~\cite{fj2002}, we find not only the lead-order asymptotics for the $L^2$-distance $\|Y_n - Y\|_2$, but in fact an exact expression for general~$n$.  It is interesting to note that the rate $n^{ - 1 / 2} (\log n)^{1 / 2}$ for $L^2$-convergence is larger even than the upper-bound rate of $n^{ - 1 / 2}$ for $d_2$-convergence from \refT{thfjd2}.

\begt[{\bf main theorem}]
\label{L2ratek} 
For $n \geq 0$ we have
$$
\|Y_n - Y\|_2^2 = (n+1)^{-1}\left(2H_n + 1 + \frac{6}{n+1}\right)  - 4\sum_{k=n+1}^\infty k^{-2} 
= 2 \frac{\ln n}{n} + O\left(\frac{1}{n}\right),
$$
where $H_n := \sum_{j=1}^n j^{-1}$ is the $n$\emph{th} harmonic number and the asymptotic expression holds as 
$n \to \infty$.
\ent

The remainder of this extended abstract is devoted to a proof of \refT{L2ratek},
which is completed in \refS{closedan}.

\section{Preliminaries}\label{pprelL2rtk}

In this section we provide recursive representations of $Y_n$ (for general~$n$) and~$Y$ that will be useful in proving \refT{L2ratek}.  Our first proposition concerns the limit~$Y$ and gives a sample-pointwise extension of the very well known~\cite{r1991} distributional identity satisfied by~$Y$.  Recall the notation~\eqref{gtheta} and the definition of~$Y$ in \refT{throsl} as the infinite series $\sum_{j = 0}^{\infty} \sum_{|\theta| = j} G_{\theta}$ in $L^2$.

\begp
\label{recYk} 
There exist random variables $F_\theta$ and $H_\theta$ for $\theta \in \Theta$ such that
\begin{enumerate}
\item the joint distributions of $(G_\theta:\theta \in \Theta)$, of $(F_\theta:\theta \in \Theta)$, and of 
$(H_\theta:\theta \in \Theta)$ agree;
 \item $(F_\theta:\theta \in \Theta)$ and $(H_\theta:\theta \in \Theta)$ are independent;
\item the series
\begin{equation}
\label{Y01}
Y^{(0)} := \sum_{j=0}^\infty\sum_{|\theta| = j}F_\theta \qquad \mbox{\rm and} \qquad 
Y^{(1)}:= \sum_{j=0}^\infty\sum_{|\theta| = j}H_\theta
\end{equation}
converge in $L^2$;
\item the random variables $Y^{(0)}$ and $Y^{(1)}$ are independent, each with the same distribution as~$Y$, and
\begin{equation}
\label{Yid}
Y  =  C(U) + U Y^{(0)} + \overline{U} Y^{(1)}.
\end{equation}
Here $U := U_1$, with $\overline{U} := 1 - U_1$, and~$C$ is defined at~\eqref{Cdef}. 
\end{enumerate}
\enp

\begin{proof}
Recall from~\eqref{gtheta} that
$$
G_{\theta} = \phi_{\theta} - 2\phi_{\theta} \ln \phi_{\theta} + 2\phi_{\theta 0} \ln \phi_{\theta 0} 
+ 2\phi_{\theta 1} \ln \phi_{\theta 1}.
$$
For $\theta \in \Theta$, define the random variable $\varphi_\theta$ (respectively, $\psi_\theta$) by
$$
\varphi_{\theta} := \phi_{0 \theta} / U \qquad (\text{resp.,}\,\, \psi_{\theta} := \phi_{1 \theta} /  \overline{U}).
$$
Then~$U$ and $\varphi_\theta$ are independent (resp.,\ $\overline{U}$ and $\psi_\theta$ are independent), 
$\varphi_{\theta}$ and $\psi_{\theta}$ each have the same distribution as $\phi_{\theta}$, and
$$
G_{0 \theta} = U F_{\theta} \qquad \mbox{and} \qquad G_{1 \theta} = \overline{U} H_{\theta},
$$
where 
\begin{align*}
F_\theta &:= \varphi_{\theta} - 2\varphi_{\theta} \ln \varphi_{\theta} + 2\varphi_{\theta 0} \ln \varphi_{\theta 0} 
+ 2\varphi_{\theta 1} \ln \varphi_{\theta 1}, \\ 
H_\theta &:= \psi_{\theta} - 2\psi_{\theta} \ln \psi_{\theta} + 2\psi_{\theta 0} \ln \psi_{\theta 0} 
+ 2\psi_{\theta 1} \ln \psi_{\theta 1}.
\end{align*}

The proposition follows easily from the clear equality
$$
\Lc(F_\theta:\theta \in \Theta) = \Lc(G_\theta:\theta \in \Theta) = \Lc(H_\theta:\theta \in \Theta),
$$
of joint laws and the fact that
\begin{eqnarray*}
Y & = & \sum_{j=0}^\infty\sum_{|\theta| = j}G_\theta = G_\varepsilon + \sum_{j=0}^\infty\sum_{|\theta| = j}G_{0 \theta} + \sum_{j=0}^\infty\sum_{|\theta| = j}G_{1 \theta}\\
& = & C(U) + U \sum_{j=0}^\infty\sum_{|\theta| = j}F_\theta + \overline{U} \sum_{j=0}^\infty\sum_{|\theta| = j}H_\theta\\
& = & C(U) + U Y^{(0)} + \overline{U} Y^{(1)}.
\end{eqnarray*}
\end{proof}

We next proceed to provide an analogue [namely, \eqref{Ynid}] of~\eqref{Yid} for each $Y_n$, rather than~$Y$, but first we need a little more notation.

Given $0 \le x < y \le 1$, let $(U_n^{xy})_{n\ge 1}$ be the subsequence of $(U_n)_{n\ge 1}$ that falls in $(x,y)$. The random variable $K_n(x,y)$ is defined to be the (random) number of key comparisons used to sort 
$U_1^{xy},\ldots,U_n^{xy}$ using \QuickSort.  The distribution of $K_n(x,y)$ of course does not depend on $(x, y)$.

We now define the random variable
\begin{eqnarray}
\label{tsndef2}
Y_{n,\theta} & := & [K_{\nu_\theta(n)}(L_\theta, R_\theta) - \mu_{\nu_\theta(n)}] / [\nu_{\theta}(n) + 1],
\end{eqnarray} 
with the centering here motivated by the fact that $\mu_{\nu_\theta(n)}$ is the conditional expectation of 
$K_{\nu_\theta(n)}(L_\theta, R_\theta)$ given $(\nu_\theta(n),L_\theta, R_\theta)$.  Then for $n \geq 1$ we have
\begin{equation}
\label{Ynid}
Y_n = \frac{n}{n + 1}\,C_n(\nu_0(n) + 1) + \frac{\nu_0(n) + 1}{n+1}\,Y_{n, 0} + \frac{\nu_1(n) + 1}{n+1}\,Y_{n, 1},
\end{equation}
where, as in~\cite{Quickapp2004}, for $1 \leq i \leq n$ we define
\begin{equation*}
C_n(i) := \mbox{$\frac{1}{n}$}(n - 1 + \mu_{i-1}+\mu_{n - i}-\mu_n).
\end{equation*}
We note for future reference that the classical divide-and-conquer recurrence for $\mu_n$ asserts precisely that
\begin{equation}
\label{sumcni0}
\sum_{i=1}^n C_n(i) = 0
\end{equation}
for $n \geq 1$.

It follows from~\eqref{Yid} and~\eqref{Ynid} that for $n \geq 1$ we have
\begin{align}
Y_n - Y
&= \left[\frac{\nu_0(n) + 1}{n+1}\,Y_{n, 0} -U Y^{(0)} \right] + \left[\frac{\nu_1(n) + 1}{n+1}\,Y_{n, 1} 
- \overline{U}Y^{(1)} \right] \nonumber \\
&{} \qquad \qquad + \left[\frac{n}{n + 1}\,C_n(\nu_0(n) + 1) - C(U)\right]\nonumber\\
  & =: W_1 + W_2 + W_3.
\label{W's}
\end{align}
Conditionally given~$U$ and $\nu_0(n)$, the random variables $W_1$ and $W_2$ are independent, each with vanishing mean, and $W_3$ is constant. Hence
$$
\EE[(Y_n - Y)^2\,\vert\,U,\nu_0(n)] = \EE[W_1^2\,\vert\,U,\nu_0(n)] + \EE[W_2^2\,\vert\,U,\nu_0(n)] + W_3^2
$$
and thus, taking expectations and using symmetry, for $n \geq 1$ we have
\begin{equation}
\label{an2decomp}
a_n^2 :=  \EE(Y_n - Y)^2  =  \EE\,W_1^2 + \EE\,W_2^2 + \EE\,W_3^2
 =  2 \EE\,W_1^2 + \EE\,W_3^2.
\end{equation}
Note that 
\begin{equation}
\label{sigmadef}
a_0^2 = \EE\,Y^2 = \sigma^2 := \mbox{$7 - \frac{2}{3}\pi^2$}
\end{equation}
(for example, \cite{Quickapp2004}).

\section{Analysis of $\EE\,W_1^2$}
\label{EW12}
In this section we analyze $\EE\,W_1^2$, producing the following result.
Recall the definition of $\sigma^2$ at~\eqref{sigmadef}.

\begp
\label{anthmk} 
Let $n \geq 1$.  For $W_1$ defined as at \eqref{W's}, we have
\begin{eqnarray*}
\EE\,W_1^2 = \frac{1}{n(n+1)^2}\sum_{k=0}^{n-1}(k+1)^2a_k^2 + \frac{\sigma^2}{6(n+1)}.
\end{eqnarray*}
\enp

For that, we first prove the following two lemmas.

\begl
\label{anthmk1} 
For any $n \ge 1$, we have
\begin{eqnarray*}
\EE \left[\left(\frac{\nu_0(n) + 1}{n+1}\right)^2\left(Y_{n, 0} - Y^{(0)} \right)^2\right] & = & \frac{1}{n}\sum_{k=0}^{n-1}\left(\frac{k+1}{n+1}\right)^2a_k^2.
\end{eqnarray*}
\enl

\begl
\label{anthmk2} 
For any $n \ge 1$, we have
\begin{eqnarray*}
\EE \left[\left(\frac{\nu_0(n) + 1}{n+1} - U\right)^2 \left( Y^{(0)} \right)^2\right] & = & \frac{\sigma^2}{6(n+1)}.
\end{eqnarray*}
\enl

\begin{proof}[\textbf{Proof of \refL{anthmk1}}] There is a probabilistic copy $Y^* = (Y^*_n)$ of the stochastic process~$(Y_n)$ such that
$$
Y_{n, 0} \equiv Y^*_{\nu_0(n)}
$$
and $Y^*$ and $Y^{(0)}$ are independent of $(U, \nu_0(n))$. This implies
\begin{eqnarray*}
\EE \left[\left(\frac{\nu_0(n) + 1}{n+1}\right)^2\left(Y_{n,0} - Y^{(0)} \right)^2\right]
 & = & \EE\,\left[\left(\frac{\nu_0(n) + 1}{n+1}\right)^2\left(Y^*_{\nu_0(n)} - Y^{(0)} \right)^2\right].
\end{eqnarray*}
By conditioning on $\nu_0(n)$, which is uniformly distributed on $\{0, \dots, n - 1\}$, we get
\begin{align*}
\EE \left[\left(\frac{\nu_0(n) + 1}{n+1}\right)^2\left(Y_{n, 0} - Y^{(0)} \right)^2\right]
 &= \EE \left[\left(\frac{\nu_0(n) + 1}{n+1}\right)^2a_{\nu_0(n)}^2\right]\\
 &= \frac{1}{n}\sum_{k=0}^{n-1}\left(\frac{k+1}{n+1}\right)^2a_k^2.
\end{align*}
\end{proof}

\begin{proof}[\textbf{Proof of \refL{anthmk2}}]
Conditionally given $\nu_0(n)$ and $Y^{(0)}$, we have that $U$ is distributed as the order statistic of rank 
$\nu_0(n) + 1$ from a sample of size~$n$ from the uniform$(0, 1)$ distribution, namely, 
Beta$(\nu_0(n)+1,n-\nu_0(n))$, with expectation $[\nu_0(n) + 1] / (n + 1)$ and variance 
$[(\nu_0(n) + 1)(n-\nu_0(n))] / [(n+1)^2(n+2)]$.  So, using also the independence of $\nu_0(n)$ and $Y^{(0)}$, we find
\begin{align*}
\lefteqn{\EE \left[\left(\frac{\nu_0(n) + 1}{n+1} - U\right)^2 \left( Y^{(0)}\right)^2\right] } \\
&= \EE \left[\frac{(\nu_0(n) + 1)(n-\nu_0(n))}{(n+1)^2(n+2)} \left( Y^{(0)} \right)^2\right]
= \sigma^2\,\EE \left[\frac{(\nu_0(n) + 1)(n-\nu_0(n))}{(n+1)^2(n+2)}\right]\\
&= \frac{\sigma^2}{(n+1)^2(n+2)}\times\frac{1}{n}\sum_{k=0}^{n-1}(k+1)(n-k)
= \frac{\sigma^2}{n(n+1)^2(n+2)}\times\frac{1}{6}n(n+1)(n+2) \\
&= \frac{\sigma^2}{6(n+1)}.
\end{align*}
\end{proof}

\begin{proof}[\textbf{Proof of \refP{anthmk}}] 
We have
\begin{align*}
\EE\,W_1^2 
&= \EE \left[\frac{\nu_0(n) + 1}{n+1}\,Y_{n, 0} - U Y^{(0)}\right]^2 \\
&= \EE \left[\frac{\nu_0(n) + 1}{n+1}\left(Y_{n, 0} - Y^{(0)} \right) + \left(\frac{\nu_0(n) + 1}{n+1} - U\right)Y^{(0)}\right]^2 \\
&= \EE \left[\left(\frac{\nu_0(n) + 1}{n+1}\right)^2\left(Y_{n, 0} - Y^{(0)} \right)^2\right] \\
&{} \qquad + \EE \left[\left(\frac{\nu_0(n) + 1}{n+1} - U\right)^2(Y^{(0)})^2\right]\\
&{} \qquad + 2\,\EE \left[\frac{\nu_0(n) + 1}{n+1}\left(Y_{n, 0} - Y^{(0)}\right)\left(\frac{\nu_0(n) + 1}{n+1} - U\right)Y^{(0)}\right].
\end{align*}
The result follows from Lemmas \ref{anthmk1}--\ref{anthmk2}, and the fact that, conditionally given $(\nu_0(n), Y_{n, 0}, Y^{(0)})$, the random variable $U$ is distributed Beta$(\nu_0(n)+1,n-\nu_0(n))$, so that the last expectation in the preceding equation vanishes.
\end{proof}

\section{Analysis of $\EE\,W_3^2$}\label{EW32}
In this section we analyze $\EE\,W_3^2$, producing the following result.

\begp
\label{W3thmk} 
For any $n \ge 1$ we have
$$
b_n^2 :=  \EE\,W_3^2  =  \frac{\sigma^2 - 7}{3} + \frac{4}{3} \left(\frac{n+2}{n+1}\right) H_n^{(2)} 
+ \frac{4}{3n(n+1)^2} H_n,
$$
where $H_n = \sum_{j=1}^n j^{-1}$ is the $n$\emph{th} harmonic number and
$H_n^{(2)}:= \sum_{j=1}^n j^{-2}$ is the $n$\emph{th} harmonic number of the second order.
\enp

For that, we first prove the following two lemmas.

\begl\label{betalog}
For any $1\leq k \leq n$ we have
$$
D(n,k) := \frac{1}{B(k,n-k+1)}\int^1_0\!t^{k-1}(1-t)^{n-k} (\ln t)\,dt = H_{k-1} - H_n,
$$
where~$B$ is the beta function.
\enl

\begl\label{cnict} For any $n \geq 1$ we have
$$
\EE[C_n(\nu_0(n)+1)C(U)] = \frac{n}{n+1} \EE[C_n(\nu_0(n)+1)]^2.
$$
\enl

\begin{proof}[\textbf{Proof of Lemma \ref{betalog}}]
The result can be proved for each fixed $n \geq 1$ by backwards induction on~$k$ and integration by parts, but we give a simpler proof.  Recall the defining expression
$$
B(\ga, \beta) = \int^1_0\!t^{\ga - 1} (1 - t)^{\beta - 1}\,dt
$$
for the beta function when $\alpha, \beta > 0$.  Differentiating with respect to~$\alpha$ gives
$$
\int^1_0\!t^{\ga - 1} (1 - t)^{\beta - 1} (\ln t)\,dt = B(\ga, \beta) [\psi(\ga) - \psi(\ga + \beta)],
$$
where~$\psi$ is the classical digamma function, i.e.,\ the logarithmic derivative of the gamma function.
But it is well known that $\psi(j) = H_{j - 1}$ for positive integers~$j$, so the lemma follows by setting
$\ga = k$ and $\beta = n - k + 1$.
%
%
\end{proof}

\begin{proof}[\textbf{Proof of Lemma \ref{cnict}}]
We know that $\nu_0(n)+1 \sim \text{unif}\{1,2,\ldots,n\}$ and that, conditionally given $\nu_0(n)$, 
the random variable~$U$ has the
Beta$(\nu_0(n)+1,n-\nu_0(n))$ distribution. So from \refL{betalog}, repeated use of~\eqref{sumcni0}, and the very well known and easily derived explicit expression
$$
\mu_n =2 (n + 1) H_n - 4 n, \qquad n \geq 0,
$$
we have
\begin{align*}
\lefteqn{\EE[C_n(\nu_0(n)+1) C(U)]} \\ 
&= \frac{1}{n} \sum_{j=1}^n C_n(j) \frac{1}{B(j,n-j+1)} \int^1_0\!t^{j-1}(1-t)^{n-j} C(t)\,dt \\
&= \frac{1}{n} \sum_{j=1}^n C_n(j) [1 + 2 \frac{j}{n+1} (H_j -H_{n+1}) + 2 \frac{n-j+1}{n+1} (H_{n-j+1} -H_{n+1})] \\
&= \frac{1}{n(n+1)} \sum_{j=1}^n C_n(j) [2 j H_j + 2 (n-j+1) H_{n-j+1}] \\
&= \frac{1}{n(n+1)} \sum_{j=1}^n C_n(j) [2jH_{j-1} - 4(j-1) + 2 (n-j+1) H_{n-j} -4 (n - j)] \\
&= \frac{1}{n(n+1)} \sum_{j=1}^n C_n(j) [\mu_{j-1} + \mu_{n-j}] 
 = \frac{1}{n(n+1)} \sum_{j=1}^n C_n(j) [\mu_{j-1} + \mu_{n-j} - \mu_n] \\
&= \frac{n}{n+1} \times \frac{1}{n} \sum_{j=1}^n C_n(j)^2 = \frac{n}{n+1} \EE[C_n(\nu_0(n)+1)]^2,
\end{align*}
as desired.
\end{proof}

\begin{proof}[\textbf{Proof of \refP{W3thmk}}] 
It follows from \refL{cnict} that
\begin{align*}
b_n^2 &= \EE \left[\frac{n}{n+1}C_n(\nu_0(n)+1) - C(U)\right]^2\\
&= \left(\frac{n}{n+1}\right)^2 \EE\,[C_n(\nu_0(n)+1)]^2 -2 \left(\frac{n}{n+1}\right) \EE[C_n(\nu_0(n)+1)C(U)] 
+ \EE\,C(U)^2\\
&= \EE\,C(U)^2 - \left(\frac{n}{n+1}\right)^2 \EE\,[C_n(\nu_0(n)+1)]^2.
\end{align*}
Knowing that $\EE\,C(U)^2 = \sigma^2 / 3$, and from the proof of Lemma A.$5$ in \cite{Quickapp2004} that
\begin{align*}
\EE\,[C_n(\nu_0(n)+1)]^2 
&= \mbox{$\frac{7}{3}\left(1 + \frac{1}{n}\right)^2 - \frac{4}{3}\left(1 + \frac{2}{n}\right)\left(1 + \frac{1}{n}\right)H^{(2)}_n - \frac{4}{3}n^{-3}H_n$},
\end{align*}
we have
$$
b_n^2 = \frac{\sigma^2 - 7}{3} + \frac{4}{3} \left( \frac{n+2}{n+1} \right) H_n^{(2)} + \frac{4}{3 n (n+1)^2} H_n,
$$
as claimed.
\end{proof}

\section{A closed form for $a_n^2$}
\label{closedan}

In this final section we complete the proof of \refT{L2ratek}, for which we need one more lemma.

\begl
\label{harm2} 
For $H_n^{(2)} = \sum_{j=1}^n j^{-2}$, the $n$\emph{th} harmonic number of the second order, we have
\begin{eqnarray*}
\sum_{j=1}^nH_j^{(2)} & = & (n+1)H_n^{(2)} - H_n
\end{eqnarray*}
for any nonnegative integer $n$.
\enl
\noindent
The lemma is well known and easily proved.
%
\smallskip

\begin{proof}[\textbf{Proof of main Theorem \ref{L2ratek}}] 
For $n \geq 1$ we have from the decomposition~\eqref{an2decomp} and Propositions~\ref{anthmk} and~\ref{W3thmk} that
\begin{align*}
a_n^2 
&= \frac{2}{n(n + 1)^2} \sum_{k = 0}^{n - 1} (k + 1)^2 a_k^2 \\
&{} \qquad + \frac{\sigma^2}{3} \left( \frac{n + 2}{n + 1} \right) -\frac{7}{3} + \frac{4}{3} \left( \frac{n + 2}{n + 1} \right)
H_n^{(2)} + \frac{4}{3 n (n + 1)^2} H_n,
\end{align*}
and we recall from~\eqref{sigmadef} that $a_0^2 = \sigma^2$.
Setting $x_n := (n + 1)^2 a_n^2$, we have $x_0 = \sigma^2$ and
$$
x_n = \frac{2}{n} \sum_{k = 0}^{n - 1} x_k + c_n \qquad \mbox{for $n \geq 1$},
$$
with 
$$
c_n := \frac{\sigma^2}{3}(n + 2)(n + 1) -\frac{7}{3}(n + 1)^2 + \frac{4}{3} (n + 2) (n + 1) H_n^{(2)} + \frac{4}{3 n} H_n.
$$
This is a standard divide-and-conquer recurrence relation for $x_n$, with solution
$$
x_n = (n + 1) \left[ \sigma^2 + \sum_{k = 1}^n \frac{k c_k - (k - 1) c_{k - 1}}{k (k + 1)} \right], \qquad n \geq 0.
$$
After straightforward computation involving the identity in \refL{harm2}, one finds
\begin{align*}
a_n^2 
&=  (n + 1)^{-1} \left( 2H_n + 1 + \frac{6}{n + 1} \right) + \sigma^2 - 7 + 4 H_n^{(2)} \\
&= (n+1)^{-1}\left(2H_n + 1 + \frac{6}{n+1}\right)  - 4\sum_{k=n+1}^\infty k^{-2} 
= 2 \frac{\ln n}{n} + O\left(\frac{1}{n}\right),
\end{align*}
as claimed.
\end{proof}

\bibliographystyle{plain}
\bibliography{references}
\end{document}